\documentclass[11pt]{article}
\usepackage{a4wide}
\usepackage{amsmath}

\usepackage{amsfonts}
\usepackage{amssymb}
\usepackage{euscript}
\newenvironment{proof}{\noindent {\it Proof.~~}\ }{\  \rule{1mm}{2mm}\medskip}

\newcommand{\A}{{\mathcal A}}
\newenvironment{proof*}{\noindent {\it Proof.~~}\ }{}

\newtheorem{theorem}{Theorem}
\newtheorem{lemma}[theorem]{Lemma}
\newtheorem{corollary}[theorem]{Corollary}
\newtheorem{proposition}[theorem]{Proposition}
\newtheorem{conj}[theorem]{Conjecture}

\begin{document}
\title{On the size of spheres of relations with a transitive group of automorphisms}
\author{Yahya Ould Hamidoune\thanks{Universit\'e Pierre et Marie Curie, Paris.
{\tt yha@ccr.jussieu.fr}}}
\maketitle

\begin{abstract} 

Let $\Gamma =(V,E)$  be a point-transitive reflexive  relation. Let $v\in V$ and   
put $r=|\Gamma (v)|.$
Also assume $\Gamma ^j(v)\cap \Gamma ^{-}(v)=\{v\}$. Then 
 $$
|\Gamma ^{j}  (v)\setminus \Gamma ^{j-1}  (v)| \ge r-1.$$

In particular we have $
|\Gamma ^{j}  (v)| \ge 1+( r-1)j.$ The last result confirms a recent conjecture of Seymour
in the case vertex-transitive graphs. Also it gives a short proof for the validity of the Caccetta-H\"aggkvist conjecture for vertex-transitive graphs and generalizes an additive result of Shepherdson.  
 
 \end{abstract}

\section{Introduction}
  Let $G$ be an abelian group  and let $S$ be a subset of $ G\setminus \{0\}$. Improving results by Chowla \cite{Chowla}, Shepherdson proved in 1947 that  
there are $s_1, \cdots , s_k\in S$ such that $k\leq \lceil \frac{|G|}{|S|}\rceil $ and  $\sum \limits_{1\leq i \leq k} s_i=0$.

In 1970 Behzad, Chartrand and  Curtis \cite{Behzad} independently conjectured that the order of $r$-biregular ($r$-outregular and inregular) directed graph $D$ without loops is $\ge r(g-1)+1$, where $g$ denotes the girth of $D$ (the smallest directed cycle in $D$).
In 1978,  Caccetta and H\"aggkvist made the stronger conjecture that the order of $r$-outregular directed graph $D$ without loops is $\ge r(g-1)+1$, where $g$ denotes the girth of $D$. These conjectures are still largely open, even for the special case $g=4$. The reader may find references and  results about this question in \cite{Bondy}. These conjectures are proved by the author in 1981 for vertex-transitive graphs \cite{HEJC}. This result applied to Cayley graphs shows the validity of Shepherdson's Theorem for all finite groups. Unfortunately we were not aware at that moment of Shepherdson's result.
More recently  Seymour proposed the following conjecture \cite{Seymour}:

Let $D$ be a directed  graph without loops and with girth $=g$. Then there is a vertex $a$ such that  $$\Gamma (a)\cup \Gamma ^2 (a)\cup \cdots \cup \Gamma ^{g-2}(a) \ge r(g-2).$$
The case   $g=4$ of this conjecture is reported in \cite{Bondy}.
 Seymour's Conjecture implies the conjecture Behzad, Chartrand and  Curtis.  One of the usual formulations of the Caccetta-H\"aggkvist Conjecture is the following : An $r$-outregular directed graph $D=(V,E)$ without loops contains a directed cycle with cardinality $\leq \lceil \frac{|V|}{r} \rceil$. Seymour's Conjecture  also implies that 
an  $r$-outregular directed graph $D=(V,E)$ without loops contains a directed cycle with cardinality $\leq \lceil \frac{|V|-1}{r} \rceil+1$.

For some technical reasons we shall use  loops. This convention is unusual in this part of Graph Theory. So we shall work with relations. Our terminology will be developed in the next section.
The diagonal of $V$ will be denoted by  $\Delta _V$.
Seymour's Conjecture may be stated in relations  language as follows:

\begin{conj} \cite{Seymour} Let $\Gamma =(V,E)$ be a reflexive  relation. 
Let $g$ be the girth of $\Gamma \setminus \Delta _V$.
 There is $x\in V$ such that  $|\Gamma ^{g-2}(x)|\ge 1+(g-2)(|\Gamma (x)|-1)$.
\end{conj}

Our main result is the following one:

Let $\Gamma =(V,E)$  be a point-transitive reflexive  relation. Let $v\in V$ and   
put $r=|\Gamma (v)|.$
Also assume $\Gamma ^j(v)\cap \Gamma ^{-}(v)=\{v\}$. Then 
 $$
|\Gamma ^{j}  (v)\setminus \Gamma ^{j-1}  (v)| \ge r-1.$$

This result implies the validity of the above conjectures for vertex-transitive graphs.

\section{Terminology }
Let $E \subset V\times V$.  The ordered pair $\Gamma = (V,E)$ will be called  a  {\em relation }. 
The relation $\Gamma$ is said to be {\em
reflexive} if $\{(x,x) \  \Big| \ x\in V
\}\subset E.$ 

Let $ a\in V$ and let $A\subset V$. We shall write $$\Gamma (a)=\{y : \ \mbox{ there is } \ x \ \mbox{ with  }\ (a,x)\in E\}$$ and $$\Gamma (A)=\bigcup \limits_{x\in A} 
\Gamma (x).$$ The cardinality of the image of $x$ will be call the degree of $x$. We shall write $d(x)=|\Gamma (x)|$. The relation $\Gamma$ will be called {\em regular} with degree $r$ if 
the elements of $V$ have  degree $=r$. 
The {\em reverse   } relation of $\Gamma $ is by definition $\Gamma ^-=(V,E^-)$, where 
$E^-=\{(x,y) \Big| \   (y,x) \in  E\}.$ 
The {\em restriction} of $\Gamma =(V,E)$ to a subset $W\subset V$ is defined as the relation 
$\Gamma _W=(W,E\cap (W\times W))$. 

Let $\Phi=(W,F)$ be a relation. A function $f : V  \longrightarrow W$ will be called a {\em homomorphism }if for all $x,y\in V$ such that $(x,y)\in E$, we have 
$(f(x),f(y))\in F$. 
The group of automorphisms of $\Gamma $ will be denoted $\A (\Gamma)$.
We shall say that a subgroup $G$ of $\A (\Gamma)$ {\em acts transitively } of $V$ if for all $x,y\in V$, there is an automorphism $f\in G$ 
such that $y=f(x)$.
The relation $\Gamma$ will be called {\em point-transitive} if $\A (\Gamma)$ acts transitively on $V$.
Clearly a point-transitive relation is regular.

One may identify a relation with its graph. In this case we mention some differences between 
our terminology (which follows closely the standard notations of Set Theory) and the notations used in some text books of Graph Theory.  We point out that 
our graphs are usually 
called directed graphs without multiple arcs. 
Notice that the of image $\Gamma (a)$ used here and in Set theory is written $\Gamma ^+(a)$ in 
some text books in Graph theory. Also our notion of degree is called outdegree. We made the choice of 
Set Theory terminology since some parts of this paper could have some interest in Group Theory and Number Theory.

We shall use the composition  $\Gamma _1 \circ\cdots \circ \Gamma _k$ of relations $\Gamma _1, 
\cdots , \Gamma _k$ on $V$. If all these relations are equal to $\Gamma$, we shall write

$$\Gamma _1 \circ \cdots \circ \Gamma _k=\Gamma ^k.$$

We shall write $\Gamma ^{0}$ for the identity relation. Also we shall write $\Gamma ^{-j}$ instead of 
$({\Gamma ^{-}})^j.$

 Let $\Gamma = (V,E)$ be relation. The girth of $\Gamma$ is by definition 

$$g(\Gamma)=\min \{k : \mbox{ there is }a \ \mbox{with} \ a\in \Gamma ^k (a)\}, $$
where $\min \emptyset =\infty$.
Notice that $g(\Gamma )$ represents the minimum size of directed cycle of the graph $\Gamma$.

We write $\Delta _V =\{(x,x) : x\in V\}$.

\section{Connectivity}

Let $\Gamma =(V,E)$ be a  relation. We shall write  $$\partial _{\Gamma}(X)= \Gamma (X)\setminus X .$$ 
When the context is clear the reference to $\Gamma$ will be omitted.

 Let $\Gamma =(V,E)$ be a  relation. 
The {\em connectivity} of $\Gamma$ is by definition $\kappa  (\Gamma)=|V|-1$, if $E= V\times V.$ 
Otherwise 
\begin{equation}  \label{eq:kappa}
\kappa  (\Gamma)=\min  \{|\Gamma (X)\setminus X|\  \Big| \ \ X\neq \emptyset \ {\rm and }\ X\cup \Gamma (X)\neq V\}.
\end{equation}

We shall say that a relation $\Gamma$ is {\em connected} if  $\kappa (\Gamma)\ge 1$.
Notice that this notion of connectedness is often called strong connectedness  in Graph Theory. Actually  we need no other
notion of connectedness in this paper. 

 For a 
relation $\Gamma$, a subset $X$ achieving the above minimum is called a
{\em fragment} of $\Gamma$. A fragment with minimal cardinality
is called an {\em atom}. The cardinality of an atom of $\Gamma$ will be denoted bu $a(\Gamma)$. It is not true that distinct  atoms are disjoint. It was proved by the author in \cite{HATOM} that either 
 distinct  atoms of $\Gamma$ are disjoint, or distinct  atoms of $\Gamma ^-$ are disjoint. As a consequence of this result we could obtain :

\begin{proposition} \cite{HATOM, HJCT} {
Let $\Gamma =(V,E)$ be a
 finite point-transitive  relation with $E\neq V\times V$ such that $ a(\Gamma)  \leq
  a(\Gamma ^-) $. 
Let $A$ be an atom of $\Gamma$. Then the subrelation $\Gamma [A]$ induced on $A$ is a 
point-transitive relation. Moreover  $|A|\le \kappa (\Gamma)$.
\label{basic}} \end{proposition}
\section{Balls}

\begin{lemma}\label{circl} Let $\Gamma =(V,E)$ be a point-transitive relation. Then for all $i$, $\Gamma ^i$ is point-transitive.  

\end{lemma}

\begin{proof}
Clearly any automorphism of $\Gamma $ is an automorphism of $\Gamma ^i$.
\end{proof}

\begin{theorem}\label{main} 
Let $\Gamma =(V,E)$  be a point-transitive reflexive  relation. Let $v\in V$ and   
put $r=|\Gamma (v)|.$
Also assume $\Gamma ^j(v)\cap \Gamma ^{-}(v)=\{v\}$. Then 
 $$
|\Gamma ^{j}  (v)\setminus \Gamma ^{j-1}  (v)| \ge r-1.$$
\end{theorem}

\begin{proof}
We shall assume  $j>1$, since  the result is obvious for $j=1.$
We may assume without loss of generality that $\Gamma$ is connected.
Clearly  $$E\neq V\times V.$$ 
 Set $\kappa =\kappa (\Gamma ).$
Let $A$ be an atom of $\Gamma$ containing   $v$.
The proof is by induction on  $d(\Gamma)=r$.

Assume first $$\kappa =r-1.$$ Then  $|\Gamma ^j(v)\setminus \Gamma ^{j-1}  (v)|\geq \kappa =r-1.$ The result holds in this case. So we may assume  
$$\kappa \leq r-2,$$ 
and hence $r\ge 3.$  Then 
$|A|\ge 2$, since otherwise 
$\kappa =|\partial (A)|=r-1$. 

{\bf Case 1}. $ a(\Gamma)  \leq
  a(\Gamma ^-) $.

  Put $K=\Gamma(A)\setminus (A\cup \Gamma (v))$.
  Put $r_*=|\Gamma (v)\cap A|$ and  $r^*=|\Gamma (v)\setminus A|$.

 By Proposition \ref{basic}, $\Gamma _A$ is point-transitive (and hence regular).
  Put $X= \Gamma ^{j-1}(v)$ and $Y=X\cup A$.

  By the induction hypothesis, we have 
  
  \begin{equation}\label{eqint}
|(\Gamma _{j}(v)\setminus \Gamma _{j-1}(v)) \cap A|\geq r_*-1. 
  \end{equation}
  
  We have clearly $\Gamma (Y)=\Gamma (X)\cup \Gamma (A\setminus X)\subset \Gamma (X)\cup \Gamma ( A\setminus v)  \subset  \Gamma (X)\cup   K$.

  It follows that $ \partial (Y)\setminus \partial (X)\subset  K.$ In particular 
  
  \begin{equation}\label{eqext}
 |\partial (Y)\cap \partial (X)|\ge |\partial (Y)|-|K|.
  \end{equation}

   Since $\Gamma ^- (v)\cap  \Gamma ^{j}(v)=\{v\},$ we  have,  $|\Gamma ^{j-1}(v)|+\kappa \le 
  |\Gamma ^{j}(v)|\le |V|-r+1$.
  It follows by Proposition \ref{basic} that $|Y|=|\Gamma ^{j-1}(v)\cup A|\leq |V|-r+1-\kappa+|A|\leq |V|-r+1.$ 

  By the definition of $\kappa$, we have $|\partial (X\cup A)|\geq \min (|V|-|Y|,\kappa)=\kappa.$

  By (\ref{eqint}) and (\ref{eqext}),  
  
  $|\partial (X)| \ge |\partial (X) \cap A|+|\partial (X) \cap \partial (Y)|\geq r_*+\kappa-|K|=r_*-1+r^*=r-1.$ 
  
{\bf Case 2}. $ a(\Gamma) >
  a(\Gamma ^-)$. The argument used in Case 1, shows that  $
|\Gamma ^{-j}  (v)\setminus \Gamma ^{-(j-1)}  (v)| \ge r-1.$
  
  By Lemma \ref{circl}, $\Gamma ^j$ is point-transitive. In particular $\Gamma ^j$
  and its reverse $\Gamma ^{-j}$ have the same degree.
  Therefore observing that these relations are reflexive 
  $$\begin{array}{lll}
   r-1\leq |\Gamma ^{-j}  (v)\setminus \Gamma ^{-(j-1)}  (v)| 
&=&|\Gamma ^{-j}  (v)|-| \Gamma ^{-(j-1)}  (v)|\\
&=&
|\Gamma ^{j}  (v)|-| \Gamma ^{(j-1)}  (v)|\\
&=&|\Gamma ^{j}  (v)\setminus \Gamma ^{(j-1)}  (v)|.
  \end{array}$$

\end{proof}

The next result shows the validity of the conjecture of Seymour mentioned in the introduction in the case of relations with a transitive group of automorphisms.

\begin{corollary}\label{Seymour} Let $\Gamma =(V,E)$  be a point-transitive reflexive   relation with 
 degree $r$ and let $v\in V$. If $\Gamma ^j(v)\cap \Gamma ^{-}(v)=\{v\}$
  then $$
|\Gamma ^{j}  (v)| \ge 1+(r-1)j.$$
\end{corollary}

\begin{proof} The proof follows by induction using Theorem \ref{main}

\end{proof}

\begin{corollary}\cite{HEJC}\label{CHVT} Let $\Gamma =(V,E)$  be a point-transitive    relation with 
 degree $r$ such that $E \cap \Delta _V=\emptyset$. Let  $g$ denotes the girth of $\Gamma $. 
  Then $
|V| \ge 1+r(g-1).$
\end{corollary}

\begin{proof} 
Set $ \Phi =(V,E \cup \Delta _V)$.
Let $v\in V$. We have clearly $\Phi ^{g-2}(v)\cap \Phi ^{-}(v)=\{v\}$.
 By Corollary  \ref{Seymour}, $|V|-r\geq |\Phi ^{g-2}(v)|\ge 1+(g-2)r.$
\end{proof}

This result, proved in \cite{HEJC}, shows the validity of Caccetta and R. H\"aggkvist Conjecture for point-transitive graphs. But the proof obtained here is much easier.

\begin{corollary}\cite{HEJC} Let $G$ be a group of order $n$ 
 and let $S\subset G\setminus \{1\}$ with cardinality $=s$. There are elements $s_1, s_2, \cdots ,
 s_k\in S$  such that $k\leq \lceil \frac{n}{s}\rceil$ and 
 $\prod \limits_{1\leq i \leq k} s_i=1$.
\end{corollary}

The proof follows by applying Corollary \ref{CHVT}, to the Cayley graph defined by $S$
on $G$. In particular the theorem of Shepherdson mentioned in the introduction holds for all finite groups.

\end{document}